\documentclass[oribibl,runningheads]{llncs} 

\usepackage{proof}
\usepackage[colorlinks=false,bookmarks=true,pdftitle={Hypersequents 
and the Proof Theory of Intuitionistic Fuzzy Logic},pdfauthor={Matthias 
Baaz and Richard Zach},pdfkeywords={many-valued logic, intuitionistic 
fuzzy logic, proof theory, Gödel-Dummett logic},dvips]{hyperref}

\def\imm#1{{\relax\ifmmode #1\else $#1$\fi}}

\spnewtheorem{thm}[theorem]{Theorem}{\bfseries}{\itshape}
\spnewtheorem{lem}[theorem]{Lemma}{\bfseries}{\itshape}
\spnewtheorem{prop}[theorem]{Proposition}{\bfseries}{\itshape}
\spnewtheorem{defn}[theorem]{Definition}{\bfseries}{\rmfamily}

\def\mms#1{{\imm{ \mathord{ 
\mathchoice{\hbox{#1}}{\hbox{#1}}{\hbox{\tiny #1}}{\hbox{\tiny 
#1}} }}}}

\def\nmodels{\not\models} 

\def\ttr{{\it tt}} 
\def\LJ{{{\bf LJ}}} 
\def\IF{{{\bf IF}}} 
\def\H{{{\bf H}}} 
\def\HIF{\imm{{\bf HIF}}} 
\def\HIFm{\imm{{\bf HIF}^-}} 
\def\HIFs{\imm{{\bf HIF}^*}}
\def\s{\mms{\bf s}}

\let\impl\supset 

\let\mmodels\models 
 
\def\LC{{\bf LC}} 

\def\len{{\rm len}}
\def\Distr{{\rm Distr}} 
\def\Frm{{\rm Frm}} 
\def\I{\imm{\Im}} 

\let\seq\vdash

\def\lorall{\imm{({\lor}{\forall})}}
\def\lorex{\imm{({\lor}{\exists})}}

\def\implex{\imm{({\impl}{\exists})}}

\makeatletter\newcommand{\ps@ref}{\addtolength{\headheight}{3ex}
\addtolength{\topmargin}{-5ex}
\addtolength{\headsep}{2ex}
\renewcommand{\@oddhead}{\footnotesize\renewcommand{\arraystretch}{.8}{
\begin{tabular}{@{}l}
Clote, Peter~G., and Helmut Schwichtenberg (eds.), \emph{Computer Science Logic.}\\
14th International Workshop, CSL 2000. Fischbachau, Germany, August 21--26, 2000.\\
Proceedings, pp.~187--201. Springer, Berlin, 2000
\end{tabular}
}} } \makeatother

\title{Hypersequents and the Proof Theory of Intuitionistic
Fuzzy Logic\protect\thanks{Research 
supported by the Austrian Science Fund under grant P--12652~MAT}} 

\titlerunning{Proof Theory of Intuitionistic Fuzzy Logic}

\author{Matthias Baaz\inst{1} \and 
\protect\href{http://www.logic.at/people/zach/}{Richard Zach}
\inst{2}}

\institute{Institut f\"ur Algebra und Computermathematik E118.2,\\
Technische Universit\"at Wien,
A--1040 Vienna, Austria,
\email{\protect\href{mailto:baaz@logic.at}{baaz@logic.at}}
\and
Institut f\"ur Computersprachen E185.2,\\
Technische Universit\"at Wien,
A--1040 Vienna, Austria,
\email{\protect\href{mailto:zach@logic.at}{zach@logic.at}}
}

\begin{document} \bibliographystyle{abbrv} 

\maketitle\thispagestyle{ref}\pagestyle{headings}

\begin{abstract}\addcontentsline{toc}{section}{Abstract}
Takeuti and Titani have introduced and investigated a logic they
called intuitionistic fuzzy logic.  This logic is characterized as the
first-order G\"odel logic based on the truth value set~$[0,1]$.  The
logic is known to be axiomatizable, but no deduction system amenable
to proof-theoretic, and hence, computational treatment, has been
known.  Such a system is presented here, based on previous work on
hypersequent calculi for propositional G\"odel logics by Avron.  It is
shown that the system is sound and complete, and allows
cut-elimination.  A question by Takano regarding the eliminability of
the Takeuti-Titani density rule is answered affirmatively.
\end{abstract}

\section{Introduction}

Intuitionistic fuzzy logic \IF{} was originally defined by Takeuti and
Titani to be the logic of the complete Heyting algebra $[0, 1]$.  In
standard many-valued terminology, \IF{} is $[0,1]$-valued first-order
G\"odel logic, with truth functions as defined below. The
finite-valued propositional versions of this logic were introduced by
G\"odel \cite{Godel:32}, and have spawned a sizeable area of logical
research subsumed under the title ``intermediate logics''
(intermediate between classical and intuitionistic logic).  The
infinite-valued propositional G\"odel logic was studied by Dummett
\cite{Dummett:59}, who showed that it is axiomatized by \LC, i.e.,
intuitionistic propositional logic plus the linearity axiom $(A \impl
B) \lor (B \impl A)$.

Takeuti and Titani \cite{TakeutiTitani:84} characterized \IF{} 
by a calculus which extends the intuitionistic predicate 
calculus \LJ{} by several axioms as well as the density rule
\[
\infer[tt']{\Gamma \seq A \lor (C \impl B)}{\Gamma \seq A \lor (C \impl p) \lor (p \impl B)} 
\]
This rule can be read as expressing the fact that the set of truth
values is densely ordered.  In this sense, the Takeuti-Titani
axiomatization is the natural axiomatization of the $[0,1]$-valued
G\"odel logic.  The valid formulas of \IF{} are also characterized as
those formulas valid in {\em every} first-order G\"odel logic based on
a linearly ordered set of truth-values (this is obvious for all logics
based on truth value sets $\subseteq [0,1]$, since a countermodel in
such a truth-value set can be straightforwardly embedded in
$[0,1]$. The general claim was established by Horn \cite{Horn:69}).
In this characterization, the density rule is not a natural
assumption, since not every linearly ordered truth-value set is densely
ordered.  It follows from this characterization that the density rule
is redundant for the axiomatization of \IF, and completeness proofs
without it have been given by Horn \cite{Horn:69} and Takano
\cite{Takano:87}.\footnote{Note that the corresponding axiom $(\forall
p)((A \impl p) \lor (p \impl B)) \impl (A \impl B)$ is not redundant
in quantified \emph{propositional} $[0,1]$-valued G\"odel logic.  See
\cite{Baaz:00a}.}  Takano posed the question of whether a
syntactic elimination of the density rule is also possible.

More recently, another axiomatizable first-order extension
of \LC{} has been studied by Corsi \cite{Corsi:89,Corsi:92} and
Avellone et al.~\cite{Avellone:99}.  This extension is defined not via
many-valued semantics but as the class of formulas valid in all
linearly ordered intuitionistic Kripke models.  It is different from
\IF; specifically, the formula \lorall{} below is not valid in it.
\IF{} can, however, also be characterized as the set of formulas valid
in all linearly ordered Kripke models with constant domains
(this was first observed by Gabbay \cite[\S3]{Gabbay:72}).

The interest of \IF{} lies in the fact that it combines properties of
logics for approximate reasoning with properties of intuitionistic
logic.  On the one hand, \IF{} is one of the basic $t$-norm logics
(see H\'ajek \cite{Hajek:98}), on the other, it is an extension of
intuitionistic logic which corresponds to concurrency (as has been
argued by Avron \cite{Avron:91b}).  We present here a calculus for
\IF{} which is adequate for further proof-theoretic study.  The basic
result in this regard is the cut-elimination theorem for this
calculus, from which a midhypersequent-theorem can be derived.  This
theorem, in turn, corresponds to Herbrand's Theorem in classical
logic, and as such is a possible basis for automated theorem proving
in~\IF.

The calculus also allows us to investigate the proof-theoretic effects
of the Takeuti-Titani rule.  We give a positive answer to Takano's
question, showing that the density rule can be eliminated from
\IF-proofs.  A simple example illustrates the possible structural
differences between proofs with and without the Takeuti-Titani rule.

\section{Syntax and Semantics of Intuitionistic Fuzzy Logic}

The language $L$ of \IF{} is a usual first-order language with
propositional variables and where free ($a$, $b$, \dots) and bound ($x$,
$y$, \dots) variables are distinguished.

\begin{defn} An {\em \IF{}-interpretation $\I = \langle D, 
\s\rangle$} is given by the {\em domain} $D$ and the {\em valuation
function} $\s$.  Let $L^D$ be $L$ extended by constants for each
element of $D$.  Then $\s$ maps atomic formulas in $\Frm(L^D)$ into
$[0,1]$, $d \in D$ to itself, $n$-ary function symbols to functions
from $D^n$ to $D$, and free variables to elements of~$D$.

The valuation function $\s$ can be extended in the obvious way to a
function on all terms.  The valuation for formulas is defined as
follows:
\begin{enumerate} 
\item $A \equiv P(t_1, \ldots, t_n)$ is atomic: 
$\I(A) = \s(P)(\s(t_1), \ldots, \s(t_n))$. \item $A \equiv \neg 
B$: \[\I(\neg B) =  \cases{0 & if $\I(B) \neq 0$ \cr                              
1 & otherwise.}\] 
\item $A \equiv B \land C$: $\I(B \land C) = \min(\I(B), \I(C))$.
\item $A \equiv B \lor C$: $\I(B \lor C) = \max(\I(A), \I(B))$.
\item $A \equiv B \impl C$: \[ \I(B \impl C) = \cases{\I(C) &
if $\I(B) > \I(C)$ \cr 1 & if $\I(B) \le \I(C)$.}\] 
\end{enumerate}
The set $\Distr_\I(A(x)) = \{\I(A(d)) : d \in D \}$ is called the {\em
distribution} of $A(x)$.  The quantifiers are, as usual, defined by
infimum and supremum of their distributions.
\begin{enumerate}
\item[(6)] $A \equiv (\forall x)B(x)$: $\I(A) = \inf 
\Distr_\I(B(x))$. \item[(7)] $A \equiv (\exists x)B(x)$: $\I(A) =  
\sup \Distr_\I(B(x))$. 
\end{enumerate} 
\I{} {\em satisfies} a formula $A$, $\I \mmodels A$, if $\I(A) = 1$. A
formula $A$ is \IF-valid if every \IF-interpretation satisfies it.
\end{defn}

Note that, as in intuitionistic logic, $\neg A$ may be defined as 
$A \impl \bot$, where $\bot$ is some formula that always takes the 
value~0.

\section{Hypersequents and IF}

Takeuti and Titani's system \IF{} is based on Gentzen's sequent
calculus \LJ{} for intuitionistic logic with a number of extra axioms
$$
\begin{array}{c}
\seq (A \impl B) \lor ((A \impl B) \impl B) \\
(A \impl B) \impl B \seq (B \impl A) \lor B \\
(A \land B) \impl C \seq (A \impl C) \lor (B \impl C) \\
(A \impl (B \lor C)) \seq (A \impl B) \lor (A \impl C)\\
(\forall x)(A(x) \lor B) \seq (\forall x)A(x) \lor B \\
(\forall x)A(x) \impl C \seq (\exists x)(A(x) \impl D) \lor (D \impl C)
\end{array}
\eqno{
\begin{array}{r}
(\textrm{Ax}1)\\
(\textrm{Ax}2)\\
(\textrm{Ax}3)\\
(\textrm{Ax}4)\\
\lorall\\
(\forall\impl)
\end{array}
}
$$
(where $x$ does not occur in $B$ or $D$) and the following additional
inference rule:
\[ 
\infer[\ttr']{\Gamma \seq A \lor (C \impl B)}{\Gamma \seq A \lor (C \impl p) \lor (p \impl B)} 
\] 
where $p$ is a propositional eigenvariable (i.e., it does not occur in
the lower sequent). It is known that the extra inference rule is
redundant.  In fact, the system \H{} of Horn \cite{Horn:69} consisting
of \LJ{} plus the schemata $$
\begin{array}{c} 
(\forall x)(A(x) \lor B) \impl (\forall x)A(x) 
\lor B \\
 (A \impl B) \lor (B \impl A) \end{array} 
\eqno{\begin{array}{r}\lorall\\
(D)\end{array}} 
$$ 
is complete for \IF{} (see also \cite{Takano:87}).  Neither of these
systems, however, has decent proof-theoretic properties such as cut
elimination, nor is a syntactic method for the elimination of the
Takeuti-Titani rule (${\it tt}'$) known.  Takano \cite{Takano:87} has
posed the question of a syntactic elimination procedure of the
Takeuti-Titani rule as an open problem.
  
We present a system which has the required properties, and which
allows the syntactic elimination of the Takeuti-Titani rule.  Our
system is based on Avron's~\cite{Avron:91b} cut-free axiomatization of
\LC{} using a hypersequent calculus.

\begin{defn} A {\em sequent} is an expression of the form \[ 
\Gamma \seq \Delta \] where $\Gamma$ and $\Delta$ are finite 
multisets of formulas, and $\Delta$ contains at most one formula. 
A {\em hypersequent} is a finite multiset of sequents, written as 
\[ \Gamma_1 \seq \Delta_1 \mid \ldots \mid \Gamma_n \seq \Delta_n 
\] 
\end{defn}

The hypersequent calculus \HIF{} has the following axioms and 
rules:

\renewcommand{\arraystretch}{1.5} 

\noindent Axioms: $A \seq A$, for any formula $A$.

\noindent Internal structural rules:
\[
\infer[iw\seq]{G \mid A, \Gamma \seq \Delta}{G \mid \Gamma \seq 
\Delta} \qquad
\infer[\seq iw]{G \mid \Gamma \seq A}{G \mid 
\Gamma \seq } \qquad 
\infer[ic\seq]{G \mid A, \Gamma \seq \Delta}
{G \mid A, A, \Gamma \seq \Delta}
\]
\noindent External structural rules:
\[\qquad
\infer[ew]{G \mid \Gamma \seq \Delta}{G}
\qquad\qquad\qquad
\infer[ec]{G \mid \Gamma \seq \Delta}
{G \mid \Gamma \seq \Delta \mid \Gamma \seq \Delta}
\]
Logical rules:
\[
\begin{array}{c@{\qquad}c} 
\infer[\neg\seq]{G \mid \neg A, \Gamma \seq}
  {G \mid \Gamma \seq A} & 
\infer[\seq\neg]{G \mid \Gamma \seq \neg A}
  {G \mid A, \Gamma \seq} \\
\infer[\lor\seq]{G \mid A \lor B, \Gamma \seq \Delta}
  {G \mid A, \Gamma \seq \Delta & G \mid B, \Gamma \seq \Delta} &
\infer[\seq\land]{G \mid \Gamma \seq A \land B}
  {G \mid \Gamma \seq A & G \mid \Gamma \seq B} \\
\infer[\seq\lor_1]{G \mid \Gamma \seq A \lor B}{G \mid \Gamma \seq A} &
\infer[\land\seq_1]{G \mid A \land B, \Gamma \seq \Delta}
  {G \mid A, \Gamma \seq \Delta} \\
\infer[\seq\lor_2]{G \mid \Gamma \seq A \lor B}{G \mid \Gamma \seq B} &
\infer[\land\seq_2]{G \mid A \land B, \Gamma \seq \Delta}
   {G \mid B, \Gamma \seq \Delta} \\
\infer[\impl\seq]{G \mid A\impl B, \Gamma_1, \Gamma_2 \seq \Delta} 
   {G \mid \Gamma_1 \seq A & G \mid B, \Gamma_2 \seq \Delta} & 
\infer[\seq\impl]{G \mid \Gamma \seq A \impl B}
   {G \mid A, \Gamma \seq B} \\
\infer[\forall\seq]{G \mid (\forall x)A(x), \Gamma \seq \Delta}
  {G \mid A(t), \Gamma \seq \Delta} & 
\infer[\seq \forall]{G \mid \Gamma \seq (\forall x)A(x)}
  {G \mid \Gamma \seq A(a)} \\
\infer[\exists\seq]{G \mid (\exists x)A(x), \Gamma \seq \Delta}
  {G \mid A(a), \Gamma \seq \Delta} & 
\infer[\seq \exists]{G \mid \Gamma \seq (\exists x)A(x)}
  {G \mid \Gamma \seq A(t)}
\end{array}
\]
Cut: \[\infer[cut]{G \mid \Gamma, \Pi \seq \Lambda}{G \mid \Gamma 
\seq A & G \mid A, \Pi \seq \Lambda}\] 
Communication: 
\[\infer[cm]{G \mid \Theta_1, \Theta_2' \seq \Xi_1 \mid 
\Theta_1', \Theta_2 \seq \Xi_2} {
G\mid \Theta_1, \Theta_1' \seq \Xi_1 & 
G\mid\Theta_2, \Theta_2' \seq \Xi_2}\] 
Density: 
\[
\infer[\ttr]{G \mid \Phi, \Psi \seq \Sigma}
{G \mid \Phi \seq p \mid p, \Psi \seq \Sigma}
\]
The rules ($\seq\forall$), $(\exists\seq)$, and $(\ttr)$ are subject
to eigenvariable conditions: the free variable $a$ and the
propositional variable $p$, respectively, must not occur in the lower
hypersequent.  We denote the calculus obtained from \HIF{} by omitting
the cut rule by \HIFm, and that obtained by omitting (\ttr) by~\HIFs.

The semantics of \IF{} can easily be extended to hypersequents
by mapping a hypersequent $H$
\[ \Gamma_1 \seq \Delta_1 \mid \ldots \mid \Gamma_n \seq \Delta_n 
\]
to the formula $H^*$
\[
(\bigwedge \Gamma_1 \impl \bigvee \Delta_1) \lor \ldots \lor 
(\bigwedge \Gamma_n \impl \bigvee \Delta_n)
\]
where $\bigwedge \Gamma_i$ denotes the conjunction of the formulas in
$\Gamma_i$ or $\top$ if $\Gamma_i$ is empty, and $\bigvee \Delta_i$
the disjunction of the formulas in $\Delta_i$ or $\bot$ if $\Delta_i$
is empty.  Deriving a formula $A$ in $\HIF$ then is equivalent to
deriving the sequent $\seq A$: the translation of $\seq A$, i.e.,
$\top \impl A$ is equivalent to $A$.

\begin{thm}[Soundness] Every hypersequent
$H$ derivable in \HIF{} is \IF-valid.
\end{thm} 

\begin{proof}
By induction on the length of the proof.  It will suffice to show that
the axioms are valid, and that the quantifier rules and (\ttr)
preserve validity.

The soundness of the quantifier rules is established by 
observing that corresponding quantifier shifting rules
are intuitionistically valid.  For instance, since
$$
\begin{array}{c}
(\exists x)(B \lor A(x)) \impl (B \lor (\exists x)A(x))\\
(\exists x)(B \impl A(x)) \impl B \impl (\exists x)A(x)\\
\end{array}
\eqno{\begin{array}{r}
\lorex\\
\implex
\end{array}}$$
are intuitionistically valid, it is easily seen that ${\seq}{\exists}$
is a sound rule.  The only problematic rules are $({\seq}{\forall})$
and $({\exists}{\seq})$.  Suppose $G \mid \Gamma \seq A(a)$ is
derivable in \HIF.  By induction hypothesis, $G^* \lor (\bigwedge
\Gamma \impl A(a))$ is valid.  Then certainly $(\forall x)(G^* \lor
(\bigwedge \Gamma \impl A(x)))$ is \IF-valid.  Since $a$ did not occur
in $G$ or $\Gamma$, we may now assume that $x$ does not either.  Since
the quantifier shift \lorall, i.e.,
\[
(\forall x)(B \lor A(x)) \impl (B \lor (\forall x)A(x)),
\]
is valid in \IF, we see that $G^* \lor (\forall x)(\bigwedge 
\Gamma \impl A(x))$ is valid. The result follows since
\[
(\forall x)(B \impl A(x)) \impl B \impl (\forall x)A(x)
\]
is intuitionistically valid, and hence \IF-valid.

The communication rule is sound as well.  Suppose the interpretation
\I{} satisfies the premises of (\textit{cm}).  The only case where
the conclusion is not obviously also satisfied is if $\I(\Theta_1')
\le \I(\Xi_1)$ and $\I(\Theta_2') \le \I(\Xi_2)$. If the left
lower sequent is not satisfied, we have $\I(\Xi_1) <
\I(\Theta_2')$, and hence $\I(\Theta_1') \le \I(\Xi_2)$, and thus
the right lower sequent is satisfied.  Similarly if the right lower
sequent is not satisfied.

For (\ttr) we may argue as follows:  Suppose that the hypersequent 
\[
H = G \mid \Phi \seq p \mid p, \Psi \seq \Sigma
\]
is \IF-valid.  Let $\I$ be an interpretation, and let $\I_r$ be just like
$\I$ except that $\I(p) = r$. Since $p$ does not occur in the
conclusion hypersequent
\[
H' = G \mid \Phi, \Psi \seq \Sigma
\]
we have $\I(H') = \I_r(H')$ and $\I(G) = \I_r(G)$.  If $\I
\models G$ we are done.  Otherwise, assume that $\I \nmodels H'$, i.e.,
\[
r_1 = \min\{\I(\Phi), \I(\Psi)\} >
\I(\Sigma) = r_2
\]
Let $r = (r_1+r_2)/2$.  Now consider $\I_r$: $\I_r \nmodels G$ by
assumption; $\I_r \nmodels \Phi \seq p$, since $\I_r(\Phi) >
r$; and $\I_r \nmodels p, \Psi \seq \Sigma$, since
$\I_r(\Psi) > r > \I_r(\Sigma)$.  Hence, $\I_r \nmodels H$, a
contradiction. \qed
\end{proof}

\begin{thm}[Completeness] Every \IF-valid hypersequent is derivable in \HIF{}.
\end{thm}

\begin{proof}
Observe that a hypersequent $H$ and its canonical translation $\seq
H^*$ are interderivable using the cut rule and the following derivable
hypersequents
\[
\begin{array}{c@{\qquad}c}
A \lor B \seq A \mid A \lor B \seq B 
& A \impl B, A \seq B \\
A \land B \seq A
& A \seq A \lor B
\end{array}
\]
Thus it suffices to show that the characteristic axioms of \IF{}
are derivable; a simple induction on the length of proofs shows
that proofs in intuitionistic predicate calculus together with
the axioms (D) and \lorall{} can be simulated in \HIF.
The formula (D) is easily derivable using the communication rule.
\[
\infer[ec]{\seq (A \impl B) \lor (B \impl A)}{
  \infer[\seq\lor]{\seq (A \impl B) \lor (B \impl A) \mid \quad 
                   \seq (A \impl B) \lor (B \impl A)}{
    \infer[\seq\lor]{\seq (A \impl B) \lor (B \impl A) \mid \quad 
                     \seq B \impl A}{
       \infer[\seq\impl]{\seq A \impl B \mid \quad\seq B \impl A}{
         \infer[\seq\impl]{\seq A \impl B \mid B \seq A}{
           \infer[\textit{cm}]{A \seq B \mid B \seq A}{
             A \seq A & B \seq B
           }
         }
       }
    }
  }
}
\]
The formula \lorall{} can be obtained thus:
{\small\[
\infer=[\seq\lor]{(\forall x)(B \lor A(x)) \seq B \lor (\forall x)A(x)}{
       \infer[\seq\forall]{(\forall x)(B \lor A(x)) \seq 
			(\forall x)A(x) \mid (\forall x)(B \lor A(x)) \seq B}{
          \infer[\forall\seq]{(\forall x)(B \lor A(x)) \seq A(a) \mid 
			(\forall x)(B \lor A(x)) \seq B}{
             \infer[\forall\seq]{(\forall x)(B \lor A(x)) \seq 
			A(a) \mid B \lor A(a) \seq B}{
                 \infer[\lor\seq]{B \lor A(a) \seq A(a) \mid 
			B \lor A(a) \seq B}{
                    \infer[\lor\seq]{B \seq A(a) \mid B \lor A(a) \seq B}{
                       \infer[cm]{B \seq A(a) \mid A(a) \seq B}{A(a) \seq A(a) & B \seq B}
                       &
                       \infer[ew]{B \seq A(a) \mid B \seq B}{B \seq B}
                    } &
                    \infer[ew]{A(a) \seq A(a) \mid B \lor A(a) \seq B}
                          {A(a) \seq A(a)}
                 }
             }
          }
       } 
} 
\]}%
The last line is obtained from the preceding by two (${\seq}{\lor}$)
inferences, followed by an external contraction.  We indicate this
with the double inference line. \qed
\end{proof}

Of course, the other axioms of Takeuti's and Titani's
system are also derivable.  We will leave the propositional axioms
1--4 as an exercise to the reader, and give the derivation on of
($\forall\impl$) as another example:
\[
\infer=[\seq\lor]{(\forall x)A(x) \impl C \seq (\exists x)(A(x) \impl D) \lor (D \impl C)}{
  \infer[\seq\impl]{\seq (\exists x)(A(x) \impl D) \mid (\forall x)A(x) \impl C \seq D \impl C}{
    \infer[\impl\seq]{\seq (\exists x)(A(x) \impl D) \mid (\forall x)A(x) \impl C, D \seq C}{
       \infer[\seq\forall]{\seq (\exists x)(A(x) \impl D) \mid D \seq (\forall x)A(x)}{
          \infer[\seq\exists]{\seq (\exists x)(A(x) \impl D) \mid D \seq A(a)}{
             \infer[\seq\impl]{\seq A(a) \impl D \mid D \seq A(a)}{
                \infer[cm]{A(a) \seq D \mid D \seq A(a)}{
                   A(a) \seq A(a) & D \seq D
                }
             }
          }
       }
       &
       \infer[ew]{\seq (\exists x)(A(x)\impl D) \mid C \seq C}{C \seq C}
    }
  }
}
\]

\section{Cut Elimination and Midhypersequent Theorem}

\begin{thm}[Cut Elimination]\label{cutel}
Any derivation of a hypersequent $G$ in \HIF{} can be transformed into
a derivation of $G$ in \HIFm{}. \end{thm} This theorem is proved in
the usual way by induction on the number of applications of the cut
rule, using the following lemma.

\begin{lem} Suppose the hypersequents 
\[
H_1 = G  \mid \Gamma \seq 
A \quad \textrm{and} \quad
H_2 = G \mid \Pi \seq \Lambda
\]
are cut-free derivable. Then \[ 
H = G \mid
\Gamma, \Pi^* \seq \Lambda
\]
where $\Pi^*$ is obtained from $\Pi$ by removing all occurrences of
$A$, is cut-free provable, and the number of applications of
(\textit{ec}) in the resulting proof is not more than the sum of
applications of (\textit{ec}) in $\gamma$ and $\delta$.
\end{lem}

\begin{proof}
Let $\gamma$ and $\delta$ be the cut-free proofs of $G$ and $H$,
respectively. We may assume, renaming variables if necessary, that
the eigenvariables in $\gamma$ and $\delta$ are distinct.  The proof
follows Gentzen's original {\em Hauptsatz}.  Define the following
measures on the pair $\langle \gamma, \delta\rangle$: the
\emph{rank}~$r = \len(\gamma) + \len(\delta)$, the \emph{degree}~$d =
\deg(A)$, and the \emph{order}~$o$ is the number of applications of
the ({\textit ec}) rule in $\gamma$, $\delta$. We proceed by induction
on the lexicographical order of $\langle d, o, r\rangle$.

If either $H_1$ or $H_2$ is an axiom, then $H$ can be derived from $H_1$ 
or $H_2$, respectively, using only weakenings. (This includes the 
case where $r = 2$).

Otherwise, we distinguish cases according to the last inferences in
$\gamma$ and $\delta$.  The induction hypothesis is that the claim of
the lemma is true whenever the degree is $< d$ or is $= d$
and either the order $< o$, or the order $=o$ and the rank $< r$.

(1) $\gamma$ or $\delta$ ends in an inference which acts on a
sequent in $G$.  We may invoke the induction
hypothesis on the premises of $H_1$ or $H_2$, and $H_2$
or $G_2$, respectively.

(2) $\gamma$ or $\delta$ ends in ($ec$).  For instance, $\gamma$ ends in
\[
\infer[ec]{G \mid \Gamma \seq A}{\infer*[\gamma']{G \mid \Gamma \seq A \mid \Gamma \seq A}{}}
\]
Apply the induction hypothesis to $\gamma'$ and $\delta$. The resulting proof $\gamma''$ of
\[
G \mid \Gamma \seq A \mid \Gamma, \Pi^* \seq \Lambda
\]
has one less ($ec$) than $\gamma$ (although it may be much longer),
and so the induction hypothesis applies again to $\gamma''$ and
$\delta$.

(3) $\gamma$ or $\delta$ end in another structural inference, (\ttr),
or (\textit{cm}): These cases are unproblematic applications of the
induction hypothesis to the premises, followed by applications of
structural inferences.

For example, assume $\gamma$ ends in (\textit{cm}), i.e.,
\[
\infer[cm]{G \mid \Theta_1, \Theta_2' \seq \Xi_1 \mid \Theta_1', \Theta_2 \seq A}{
  \infer*[\gamma_1]{G \mid \Theta_1, \Theta_1' \seq \Xi_1}{} &
  \infer*[\gamma_2]{G \mid \Theta_2, \Theta_2' \seq A}{}
}
\]
where $\Gamma = \Theta_1', \Theta_2$. Apply the deduction hypothesis to
the right premise and $H_2$ to obtain a cut-free proof of
\[
G \mid \Theta_2, \Theta_2', \Pi^* \seq \Lambda
\]
Using applications of (\textit{ew}) and (\textit{cm}),
we obtain the desired result.

The case of ($tt$) may be of special 
interest.  Suppose $\gamma$ ends in(\ttr), with 
\[
\infer[tt]{G \mid \Phi, \Psi \seq A}{G \mid \Phi \seq p \mid p, \Psi \seq A}
\]
Apply the induction
hypothesis to the premises of $H_1$ and $H_2$, and apply (\ttr) to obtain
the desired proof:
\[
\infer[tt]{G \mid \Phi, \Psi, \Pi^* \seq \Lambda}{
G \mid \Phi \seq p \mid p, \Psi, \Pi^* \seq \Lambda}
\]
The case of $\delta$ ending in ($tt$) is handled similarly.

(4) $\gamma$ ends in a logical inference not involving the cut formula,
or $\delta$ ends in a logical inference not involving the cut formula.  These
cases are easily handled by appeal to the induction hypothesis and
application of appropriate logical and structural inferences.  We
outline the case where $\gamma$ ends in $(\impl\seq)$:
\[
\infer[\impl\seq]{G \mid B \impl C, \Gamma \seq A}{
  \infer*[\gamma_1]{G \mid C, \Gamma \seq A}{} &
  \infer*[\gamma_2]{G \mid \Gamma \seq B}{}
}
\]
We apply the induction hypothesis to the left premise and $H_2$, and apply ($\impl\seq$):
\[
\infer{G \mid B \impl C, \Gamma, \Pi^* \seq \Lambda}{G \mid C, \Gamma, \Pi^* \seq \Lambda &
G \mid \Gamma \seq B}
\]

(5) Both $\gamma$ and $\delta$ end in logical inferences acting on a
cut formula.  For instance, if $A = B \impl C$ we have
\[
\infer[\seq\impl]{G \mid \Gamma \seq B \impl C}{
  \infer*[\gamma_1]{G \mid B, \Gamma \seq C}{}
}
\qquad
\infer[\impl\seq]{G \mid B \impl C, \Pi_1, \Pi_2 \seq \Lambda}{
  \infer*[\delta_1]{G \mid \Pi_1 \seq B}{} &
  \infer*[\delta_2]{G \mid C, \Pi_2 \seq \Lambda}{}
}
\]
First we find proofs $\delta_1'$ and $\delta_2'$ of
\[
G \mid \Gamma, \Pi_1^* \seq B \quad{\textrm and}\quad
G \mid C, \Gamma, \Pi_2^* \seq \Lambda
\] 
either by applying the induction hypothesis to $\gamma$ and $\delta_1$
or $\delta_2$ if $\Pi_1$ or $\Pi_2$, respectively, contain $B \impl
C$, or otherwise by adding (\textit{ic})-inferences to $\delta_1$ and
$\delta_2$.  Now apply the induction hypothesis based on the reduced
degree of the cut formulas twice: first to $\delta_1'$ and $\gamma_1$
to obtain $G \mid \Gamma, \Gamma, \Pi_1^* \seq C$, and then to the
resulting proof and $\delta_2'$ to obtain
\[
G \mid \Gamma, \Gamma, \Gamma, \Pi_1^*, \Pi_2^*
\seq \Lambda.
\] 
The desired result follows by several applications of~(\textit{ic}).

The other cases are similar and are left to the reader. \qed
\end{proof}

Cut elimination is a basic prerequisite for proof theoretic and
computational treatments of a logic.  As an immediate consequence of
cut elimination we have the subformula property: every \IF-valid
formula has a proof which only contains subformulas of the endformula
(plus possibly propositional variables used in (\ttr)).  Another
important corollary is the midhypersequent theorem.  It corresponds to
Herbrand's Theorem for classical logic and is thus the basis for any
resolution-style automated proof method.

\begin{thm}\label{midseq}
Any hypersequent~$H$ with only prefix formulas has a proof where no
propositional inference follows a quantifier inference.  Such a proof
contains one or more hypersequents~$M$, called midhypersequents, so
that $M$ contains no quantifiers, all the inferences above $M$ are
propositional or structural, and all the inferences below $M$ are
either quantifier inferences of structural inferences.
\end{thm}

\begin{proof}
This is proved exactly as for the classical and intuitionistic case
(see Takeuti \cite{Takeuti:87}).  First, observe that all axioms are
cut-free derivable from atomic axioms.  The cut-elimination theorem thus
provides us with a cut-free proof $\pi$ of $H$ from atomic axioms. Next, observe that the ($\lor\seq$) rule can be simulated without using cuts by the rule
\[
\infer[\lor\seq']{G \mid A \lor B, \Gamma \seq \Delta_1 \mid 
                         A \lor B, \Gamma \seq \Delta_2}
      {G \mid A, \Gamma \seq \Delta_1 & G \mid B, \Gamma \seq \Delta_2}
\]
The rule can be derived as follows (we omit side sequents):
\[
\infer[\lor\seq]{A \lor B, \Gamma \seq \Delta_1 \mid 
                 A \lor B, \Gamma \seq \Delta_2}{
    \infer[\lor\seq]{A \lor B, \Gamma \seq \Delta_1 \mid 
                            A, \Gamma \seq \Delta_2}{
        \infer[\textit{cm}]{B, \Gamma \seq \Delta_1 \mid 
                            A, \Gamma \seq \Delta_2}{
             A, \Gamma \seq \Delta_1 & B, \Gamma \seq \Delta_2
        }
        & A, \Gamma \seq \Delta_1
    } & B, \Gamma \seq \Delta_2
}
\]
Of course, ($\lor\seq'$) together with (\textit{ec}) simulates
($\lor\seq$).  We replace all applications of ($\lor\seq$) by
applications of ($\lor\seq'$) in our cut-free proof.

Define the order of a quantifier inference in $\pi$ to be the number
of propositional inferences under it, and the order of $\pi$ as the
sum of the orders of its quantifier inferences.  The proof is by
induction on the order of $\pi$.  The only interesting case is of
$(\lor\seq')$ occurring below a quantifier inference, since this case
does not work for intuitionistic logic.

Suppose $\pi$ contains a ($\seq\forall$) inference above a
($\lor\seq'$) inference, and so that all the inferences in between are
structural. We have the following situation:
\[
\infer[\lor\seq']{G \mid A \lor B, \Gamma \seq \Delta \mid 
                         A \lor B, \Gamma \seq (\forall x)A(x)}{
  \infer*{G \mid A, \Gamma \seq \Delta}{} &
  \infer*[\delta]{G \mid B, \Gamma \seq (\forall x)A(x)}{
     \infer[\seq\forall]{G' \mid \Gamma' \seq (\forall x)A(x)}{
        \infer*[\delta']{G' \mid \Gamma' \seq A(a)}{}
     }
  }
} 
\]
where $\delta$ contains only structural inferences.  We reduce the
order of $\pi$ by replacing this part of $\pi$ by:

\begin{center}\hfill\quad
$\infer[\seq\forall]{G \mid A \lor B, \Gamma \seq \Delta \mid 
                           A \lor B, \Gamma \seq (\forall x)A(x)}{
  \infer[\lor\seq']{G \mid A \lor B, \Gamma \seq \Delta \mid 
                           A \lor B, \Gamma \seq A(a)}{
     \infer*{G \mid A, \Gamma \seq \Delta}{} &
     \infer*[\delta]{G \mid B, \Gamma \seq A(a)}{
         \infer*[\delta']{G' \mid \Gamma' \seq A(a)}{}
     }
  }
}$
\hfill\squareforqed
\end{center}
\end{proof}

\section{Elimination of the Takeuti-Titani Rule}

The Takeuti-Titani rule is the least understood feature of the
original Takeuti-Titani axiomatization of \IF.  We show below that the
rule can be eliminated from proofs in~\HIF.  This had been posed as a
problem by Takano \cite{Takano:87}.  The proof is by induction on the
number of applications of (\ttr) and the length of the proof.  The exact
complexity of the elimination procedure is still to be investigated.
The (\ttr) rule can have significant effects on proof structure.  For
instance, one of the calculi in Avron \cite{Avron:91b} uses the split
rule
\[
\infer[\textit{split}]{G \mid \Gamma  \seq \Delta \mid 
                              \Gamma' \seq \Delta}
       {G \mid \Gamma, \Gamma' \seq \Delta}
\]
If this rule is added to \HIF, it is possible to transform proofs 
so that each application of the communication rule has a premise
which is a propositional axiom.  This is not possible without
(\ttr).  The transformation works by replacing each occurrence of
the communication rule by
{\small\[
\infer={G_1 \mid G_2 \mid \Gamma_1, \Gamma_2' \seq A_1 \mid \Gamma_1', \Gamma_2 \seq A_2}{
  \infer[tt]{G_1 \mid G_2 \mid \Gamma_1 \seq A_1 \mid \Gamma_2 \seq A_2 \mid \Gamma_2 \seq A_1 \mid \Gamma_2' \seq A_2}{
      \infer[tt]{G_1 \mid G_2 \mid \Gamma_1 \seq A_1 \mid \Gamma_2 \seq A_2 \mid p \seq A_1 \mid \Gamma_2 \seq p \mid \Gamma_2' \seq A_2}{
          \infer[cut]{G_1 \mid G_2 \mid \Gamma_1 \seq A_1 \mid \Gamma_2 \seq q \mid p \seq A_1 \mid \Gamma_2 \seq p \mid q \seq A_2 \mid \Gamma_2' \seq A_2}{
              \infer[cm]{G_1 \mid \Gamma_1 \seq A_1 \mid \Gamma_1' \seq q \mid p \seq A_1 \mid q \seq p}{
                   q \seq q &
                   \infer[cm]{G_1 \mid \Gamma_1 \seq A_1 \seq \Gamma_1' \seq p \mid p \seq A_1}{
                       p \seq p &
                       \infer[\textit{split}]{G_1 \mid \Gamma_1 \seq A_1 \mid \Gamma_1' \seq A_1}{
                           G_1 \mid \Gamma_1, \Gamma_1' \seq A_1
                       }
                   }
              } % &
              \infer[cm]{G_2 \mid \Gamma_2 \seq q \mid q \seq A_2 \mid \Gamma_2' \seq A_2}{
                  \infer[split]{G_2 \mid \Gamma_2 \seq A_2 \mid \Gamma_2' \seq A_2}{
                      G_2 \mid \Gamma_2, \Gamma_2' \seq A_2
                  } & q \seq q
              }
          }
      }
  }
}
\]}

\begin{prop}\label{noweak} Let $\delta$ be a \HIFs-derivation 
of hypersequent $H$ with length $k$, where $H$ is of the form 
\[ 
G \mid \Gamma_1, \Pi_1 \seq \Delta_1, \Pi'_1 \mid \ldots \mid \Gamma_n,
\Pi_n \seq \Delta_n, \Pi'_n 
\] 
and $\bigcup \Pi_i \subseteq \{p\}$,
$\Pi_i' = \emptyset$, and $p$ does not occur in $G$, $\Gamma_i$ or
$\Delta_i$ ($\bigcup \Pi_i' = \{p\}$, $\Pi_i = \emptyset$, and $p$
does not occur in $G$, $\Gamma_i$ or $\Delta_i$).

Then the hypersequent $G \mid \Gamma_{i_1} \seq \Delta_{i_1} \mid \ldots 
\mid \Gamma_{i_m} \seq \Delta_{i_m}$ is derivable in length $\le 
k$.\end{prop}

\begin{proof} Easy induction on $k$.  Every occurrence of $p$ must 
arise from a weakening, simply delete all these weakenings.
\end{proof}

\begin{thm} Applications of (\textit{tt}) can be eliminated from 
\HIF-derivations.
\end{thm}
This follows from the following lemma by induction on the number of
applications of (\ttr) in a given \HIFm-derivation.

\begin{lem} 
If $\delta$ is an \HIFs-derivation of 
\[ 
H = G \mid \Phi_1 \seq \Pi_1 \mid \ldots \Phi_n \seq \Pi_n \mid 
\Pi'_1, \Psi_1 \seq \Sigma_1 \mid \ldots \mid \Pi'_m, \Psi_m \seq \Sigma_m,
\] 
where $p$ does not occur in $G$, $\Phi_i$, $\Psi_i$ or $\Sigma_i$, and
$\bigcup \Pi_i \cup \bigcup \Pi'_i \subseteq \{p\}$, then there is a
\HIFs-derivation of
\[ 
H^* = G \mid 
\Phi_1, \ldots, \Phi_n, \Psi_1 \seq \Sigma_1 \mid \ldots \mid 
\Phi_1, \ldots, \Phi_n, \Psi_m \seq \Sigma_m.
\] 
\end{lem}

\begin{proof}
By induction on the length of~$\delta$.  We distinguish cases
according to the last inference~$I$ in $\delta$.  For simplicity, we
will write $p$ in what follows below instead of $\Pi_i$ or $\Pi_i'$
with the understanding that it denotes an arbitrary multiset of $p$'s.

(1) The conclusion of of $I$ is so that $p$ only occurs on the right
side of sequents, or only on the left side.  Then Prop.~\ref{noweak}
applies, and the desired hypersequent can be derived without (\ttr).

(2) $I$ applies to sequents in $G$.  Then the induction hypothesis can
be applied to the premise(s) of $I$ and appropriate inferences added
below.

(3) $I$ is structural inference other than (\textit{cut}) and
(\textit{cm}), or a logical inference with only one premise, or a
logical inference which applies to a $\Sigma_i$.  These cases are
likewise handled in an obvious manner and are unproblematic. One 
instructive example might be the case of ($\impl\seq$). Here the
premises would be of the form, say,
\[
\begin{array}{l}
G \mid \Phi_1 \seq p \mid \Phi_2 \seq p \ldots \mid \Phi_n \seq p \mid 
p, \Psi_1 \seq \Sigma_1 \mid \ldots \mid p, \Psi_m \seq \Sigma_m \mid 
p, \Gamma_1 \seq A \\
G \mid \Phi_1 \seq p \mid \Phi_2 \seq p \ldots \mid \Phi_n \seq p \mid 
p, \Psi_1 \seq \Sigma_1 \mid \ldots \mid p, \Psi_m \seq \Sigma_m \mid 
B, \Gamma_2 \seq p
\end{array}
\]
Let $\Phi = \Phi_1, \ldots, \Phi_n$. The induction hypothesis
provides us with
\[
\begin{array}{l}
G \mid \Phi, \Psi_1 \seq \Sigma_1 \mid \ldots \mid 
       \Phi, \Psi_m \seq \Sigma_m \mid \Phi, \Gamma_1 \seq A \\
G \mid B, \Gamma_2, \Phi, \Psi_1 \seq \Sigma_1 \mid \ldots \mid 
       B, \Gamma_2, \Phi, \Psi_m \seq \Sigma_m 
\end{array}
\]
We obtain the desired hypersequent by applying ($\impl\seq$) successively
$m$ times, together with some contractions.

(4) $I$ is a cut. There are several cases to consider, most of 
which are routine.  The only tricky case is when the cut formula is
$p$ and $p$ occurs both on the left and the right side of sequents
in both premises of the cut. For simplicity, let us consider the cut rule
in its multiplicate formulation 
\[
\begin{array}{l}
G \mid \Phi_1 \seq p \mid \ldots \mid \Phi_n \seq p \mid 
p, \Psi_1 \seq \Sigma_1 \mid \ldots \mid p, \Psi_m \seq \Sigma_m \mid 
\Gamma \seq p \\ 
G \mid \Phi_1 \seq p \mid \ldots \mid \Phi_n \seq p \mid 
p, \Psi_1 \seq \Sigma_1 \mid \ldots \mid p, \Psi_m \seq \Sigma_m \mid 
p, \Pi \seq \Lambda
\end{array}
\]
We want to find a derivation of
\[
G \mid \Phi, \Psi_1 \seq \Sigma_1 \mid 
\Phi, \Psi_m \seq \Sigma_m \mid \Gamma, \Pi \seq \Lambda
\]
where $\Phi = \Phi_1, \ldots, \Phi_n$.
The induction hypothesis applied to the premises of the cut gives us
\[
\begin{array}{l}
G \mid \Gamma, \Phi, \Psi_1 \seq \Sigma_1 \mid \ldots \mid 
       \Gamma, \Phi, \Psi_m \seq \Sigma_m \\ 
G \mid \Phi, \Psi_1 \seq \Sigma_1 \mid \ldots \mid 
       \Phi, \Psi_m \seq \Sigma_m \mid \Phi, \Pi \seq \Lambda
\end{array}
\]
We obtain the desired hypersequent by $m$ successive applications 
of (\textit{cm}).

(5) $I$ is ($\lor\seq$), or ($\exists\seq$) applying to
    $\Phi_i$ or $\Psi_i$. Consider the case of ($\lor\seq$), the
    others are treated similarly. The premises of $I$ are, for
    example,
\[
\begin{array}{l}
G \mid A, \Phi_1 \seq p \mid \Phi_2 \seq p \ldots \mid \Phi_n \seq p \mid 
p, \Psi_1 \seq \Sigma_1 \mid \ldots \mid p, \Psi_m \seq \Sigma_m \\ 
G \mid B, \Phi_1 \seq p \mid \Phi_2 \seq p \ldots \mid \Phi_n \seq p \mid 
p, \Psi_1 \seq \Sigma_1 \mid \ldots \mid p, \Psi_m \seq \Sigma_m 
\end{array}
\]
By induction hypothesis, we obtain
\[
\begin{array}{l}
G \mid A, \Phi_1, \ldots, \Phi_n, \Psi_1 \seq \Sigma_1 \mid \ldots \mid
 A, \Phi_1, \ldots, \Phi_n, \Psi_m \seq \Sigma_m\\
G \mid B, \Phi_1, \ldots, \Phi_n, \Psi_1 \seq \Sigma_1 \mid \ldots \mid
 B, \Phi_1, \ldots, \Phi_n, \Psi_m \seq \Sigma_m
\end{array}
\]
It is not straightforwardly possible to derive the desired
hypersequent from these.  If $\Psi_i = \{P_{i1}, \ldots, P_{ik_i}\}$,
let $Q_i = P_{i1} \impl \ldots P_{ik_i} \impl \Sigma_i$.  Then we do
easily obtain, however, the following by repeated application of
($\seq\impl$), ($\seq\lor$) and~($ec$):
\[
\begin{array}{l}
G \mid A, \Phi_1, \ldots, \Phi_n \seq Q_1 \lor \ldots \lor Q_m\\
G \mid B, \Phi_1, \ldots, \Phi_n \seq Q_1 \lor \ldots \lor Q_m
\end{array}
\]
Now a single application of ($\lor\seq$), plus (\textit{ec}) gives us 
\[
K = G \mid \underbrace{A \lor B, \Phi_1, \ldots, \Phi_n}_\Gamma 
\seq Q_1 \lor \ldots \lor Q_m
\]
Then we derive, using $m-1$ cuts:
\[
\infer={\Gamma \seq Q_1 \mid \ldots \mid \Gamma \seq Q_m}{
  \infer*{\Gamma \seq Q_1 \mid \ldots \mid \Gamma \seq Q_{m-1} \lor Q_m}{
     \infer={\Gamma \seq Q_1 \mid \Gamma \seq \underbrace{Q_2 \lor \ldots \lor Q_m}_Q}{
       K & 
       \infer*[\delta_1]{Q_1 \lor Q \seq Q_1 \mid Q_1 \lor Q \seq Q}{}
     }
   } &
   \infer*[\delta_{m-1}]{Q_{m-1} \lor Q_m \seq Q_{m-1} \mid Q_{m-1} \lor Q_m \seq Q_m}{}
}
\]
where $\delta_i$ is the derivation
\[
\infer[\lor\seq]{Q_i \lor \underbrace{Q_{i+1} \lor \ldots \lor Q_m}_Q \seq Q_i \mid Q_i \lor \ldots \lor Q_m \seq \underbrace{Q_{i+1} \lor\ldots \lor Q_m}_Q}{
  Q_i \seq Q_i &
  \infer[\lor\seq]{Q \seq Q_i \mid Q_i \lor Q \seq Q}{
     \infer[cm]{Q\seq Q_i \mid Q_i \seq Q}{
       Q \seq Q  & Q_i \seq Q_i
     } &
  Q \seq Q
  }
}
\]
The desired hypersequent is obtained by $m$ cuts with
\[
Q_i, P_{i1}, \ldots, P_{ik_i} \seq \Sigma_i
\]

(6) $I$ is a communication rule.  This is the most involved case, as
several subcases have to be distinguished according to which of the
two communicated sequents contains $p$.  Neither of these cases are
problematic.  We present two examples:

(a) One of the communicated sequents contains $p$ on the right.  Then
the premises of $I$ are
\[
\begin{array}{l}
G \mid \Phi_1 \seq p \mid \ldots \mid \Phi_n \seq p \mid 
p, \Psi_1 \seq \Sigma_1 \mid \ldots \mid p, \Psi_m \seq \Sigma_m \mid 
\Theta_1, \Theta_1' \seq p  \\
G \mid \Phi_1 \seq p\mid \ldots \mid \Phi_n \seq p \mid 
p, \Psi_1 \seq \Sigma_1 \mid \ldots \mid p, \Psi_m \seq \Sigma_m \mid 
\Theta_2, \Theta_2' \seq \Xi_2 
\end{array}
\] 
where. The induction hypothesis applies to
these two hypersequents.  If we write $\Phi = \Phi_1, \ldots,
\Phi_n$,  we have
\[
\begin{array}{l}
G\mid \fbox{$\Theta_1, \Theta_1', \Phi, \Psi_1 \seq \Sigma_1$} \mid \ldots 
\mid \Theta_1, \Theta_1', \Phi, \Psi_m \seq \Sigma_m \\
G \mid \fbox{$\Theta_2, \Theta_2' \seq \Xi$} \mid \Phi, \Psi_1 \seq \Sigma_1 
\mid \ldots \mid \Phi, \Psi_m \seq \Sigma_m \\
\end{array}
\] 
We obtain the desired result by applying $m$ instances of
(\textit{cm}), internal weakenings and external contractions as
necessary, to obtain, in sequence
\[
\begin{array}{l}
G \mid \Theta_1, \Theta_2', \Phi, \Psi_1 \seq \Sigma_1 \mid \ldots \mid 
\Theta_1, \Theta_1', \Phi, \Psi_m \seq \Sigma_m \mid 
\Theta_1', \Theta_2 \seq \Xi\\
\multicolumn{1}{c}{\ddots}\\
G \mid  \Theta_1, \Theta_2', \Phi, \Psi_1 \seq \Sigma_1 \mid \ldots \mid 
\Theta_1, \Theta_2', \Phi, \Psi_m \seq \Sigma_m \mid 
\Theta_1', \Theta_2 \seq \Xi
\end{array}
\] 
The sequents participating in the application of (\textit{cm}) are marked by
boxes.  The original end hypersequent follows from the last one by
internal weakenings.

(b) The communicated sequents both contain $p$, once on the right,
once on the left.  The premises of $I$ are
\[
\begin{array}{l}
G \mid \Phi_1 \seq p \mid \ldots \mid \Phi_n \seq p \mid 
p, \Psi_1 \seq \Sigma_1 \mid \ldots \mid p, \Psi_m \seq \Sigma_m \mid 
\Theta_1, \Theta_1' \seq p \\
G \mid \Phi_1 \seq p \mid \ldots \mid 
\Phi_n \seq p \mid p, \Psi_1 \seq \Sigma_1 \mid \ldots \mid 
p, \Psi_m \seq \Sigma_m \mid 
p, \Theta_2, \Theta_2' \seq \Xi 
\end{array}
\] 
We have proofs of
\[
\begin{array}{l}
G \mid \Theta_1, \Theta_1', \Phi, \Psi_1 \seq \Sigma_1 \mid \ldots \mid 
       \Theta_1, \Theta_1', \Phi, \Psi_m \seq \Sigma_m \\
G \mid \Phi, \Psi_1 \seq \Sigma_1 \mid \ldots \mid \Phi, \Psi_m \seq \Sigma_m 
\mid \Theta_2, \Theta_2', \Phi \seq \Xi 
\end{array}
\] 
Again, a sequence of $m$ applications of (\textit{cm}), together with
internal weakenings and external contractions produces the desired end
sequent. \qed
\end{proof}

Note that in case (5), several new cuts are introduced.  As a
consequence, the elimination procedure does not directly work for
cut-free proofs. If a proof with neither cut nor communication is
required, the elimination procedure has to be combined with the
cut-elimination procedure of Thm.~\ref{cutel}.  The additional cuts
can be avoided by replacing ($\lor\seq$) and ($\exists\seq$) by the 
following generalized rules:
\[
\begin{array}{c}
\infer[\lor\seq^*]{G \mid A \lor B, \Gamma_1 \seq \Delta_1 \mid \ldots \mid A \lor B, \Gamma_n \seq \Delta_n}{
G \mid A, \Gamma_1 \seq \Delta_1 \mid \ldots \mid A, \Gamma_n \seq \Delta_n &
G \mid B, \Gamma_1 \seq \Delta_1 \mid \ldots \mid B, \Gamma_n \seq \Delta_n}\\
\infer[\exists\seq^*]{G \mid (\exists x)A(x), \Gamma_1 \seq \Delta_1 \mid \ldots \mid (\exists x)A(x), \Gamma_n \seq \Delta_n}{G \mid A(a), \Gamma_1 \seq \Delta_1 \mid \ldots \mid A(a), \Gamma_n \seq \Delta_n}
\end{array}
\]
These rules, however, cannot be simulated by the ordinary rules
without using cut (the simulation with cut is given in case (5)).  By
changing case (5) accordingly, the elimination procedure will
transform a cut-free \HIF-derivation into a cut-free one without
(\ttr), but with ($\lor\seq^*$) and ($\exists\seq^*$).

\end{document}